\documentclass[12pt]{article}

\usepackage{graphicx}
\usepackage{amssymb}
\usepackage{amsfonts}

\usepackage{times,amsmath,amsthm,amssymb,latexsym,amscd}

\newcommand{\syst}{{\rm{syst}}}

\newcommand{\calF}{{\mathcal{F}}}

\newcommand{\eps}{\epsilon}
\newtheorem{Th}{Theorem}[section]

\newcommand{\btbt}{\left( \begin{array}{cc}}
\newcommand{\etbt}{ \end{array}\right)}
\newcommand{\bthbth}{\left( \begin{array}{ccc}}
\newcommand{\ethbth}{ \end{array}\right)}
\newtheorem{Lm}{Lemma}[section]

\newtheorem{Pro}{Proposition}[section]

\newcommand{\Mp}{\mathcal{M}_{P}}

\def\G{{\Gamma}}

\newcommand{\calO}{{\mathcal{O}}}
\newcommand{\calR}{{\mathcal{R}}}
\newcommand{\calP}{{\mathcal{P}}}
\newcommand{\calL}{{\mathcal{L}}}
\newcommand{\V}{{\cal{V}}}
\newcommand{\U}{{\cal{U}}}
\newcommand{\length}{{\hbox{length}}}
\newcommand{\calFns}{{\cal{F}}_n^*}
\newcommand{\calFn }{{\cal{F}}_n}

\newcommand{\bcol}{\left(\begin{array}{c}}
\newcommand{\ecol}{\end{array}\right)}
\renewcommand{\l}{\length}
\newcommand{\tr}{{\hbox{tr}}}

\newcommand{\Prob}{{\hbox{Prob}}}

\title{The length of closed geodesics on random Riemann Surfaces.}
\author{Eran Makover\thanks{makovere@ccsu.edu} and Jeffrey McGowan\thanks{mcgowan@ccsu.edu}}
\begin{document}
\maketitle
\begin{abstract}
Short geodesics are important in the study of the geometry and the spectra of
Riemann surfaces. Bers' theorem gives a global bound on the length of the
first $3g-3$ geodesics.
We use the construction of Brooks and Makover of random Riemann surfaces to
investigate the distribution of short ($< \log (g)$) geodesics on a
random Riemann surfaces. We calculate the expected value of the shortest
geodesic, and show that if one orders prime non-intersecting geodesics by length $\gamma_1\le \gamma_2\le \dots \le \gamma_i ,\dots$, then for fixed $k$, if one allows the genus to go to infinity, the length of $\gamma_{k}$ is independent of the genus.\end{abstract}

\section{Introduction}
A standard tool in the study of compact Riemann surfaces is the decomposition into ``pairs of pants'' (Y pieces).  Given a surface of genus $g\ge 2$, there are $3g-3$ simple closed geodesics which partition the surface into $g-1$ such pieces.  Bounds on the lengths of the geodesics in such partitions are extremely desirable.  If the geodesics are $\gamma_{1}, \ldots, \gamma_{3g-3}$, their lengths $l(\gamma_{i})$ give half of the {\it Fenchel-Nielsen} parameters which parametrize the $6g-6$ dimensional Teichm\"uller space of compact surfaces of genus $g$.

Bers (\cite{Bers1,Bers2}) proved that for every compact Riemann surface of genus $g\ge 2$ there is a partition with $$l(\gamma_{1},\ldots,l(\gamma_{3g-3})) \le L_{g}$$ where $L_{g}$ is a constant depending only on $g$.  The best possible such constant is called {\it Ber's constant}.  A constructive argument due to Abikoff (\cite{Abikoff}) gives an explicit bound for $L_{g}$; unfortunately this bound grows faster than exponentially in $g$.  The best result known is 

\begin{Th}[\cite{buser}]\label{buser}
Every compact Riemann Surface of genus $g\le 2$ has partition $\gamma_1, \dots , \gamma_{3g-3}$ satisfying 
$$l(\gamma_k) \le 4k\log \frac{8\pi(g-1)}{k} \hbox{\hskip 1in} k=1, \dots ,3g-3$$

The longest geodesic is bounded by \begin{equation}\label{bersbnd}\gamma_{3g-3} \le 26(g-1)\end{equation}
\end{Th}

One might hope that in fact the bound in (\ref{bersbnd}) might be improved to a logarithmic one, but this is impossible.  The ``hairy torus'' gives a \emph{lower bound} of 

\begin{Pro}
$L_{g}\ge \sqrt{6g} - 2$ for all $g \ge 2$.
\end{Pro}

The length of the {\it shortest} geodesic $l(\gamma_{1})$ is also of particular interest.  There are examples of classes of surfaces, such as the  conformal compactification of the principal modular surfaces, where it is known that  (\cite{buser-sarnak, platonic} ) 

$$l(\gamma_{1} )= O(\log{g}).$$

Recently  Katz  Schaps and   Vishne \cite{KSV} show for Hurwitz
surfaces  and some other principal congruence subgroups of arbitrary arithmetic surfaces a similar behavior. But all these examples are rare in the sense that they not occur  in all genera and when they do occur  they there are  a small number of such surfaces.
In this paper, we study the Belyi surfaces .  Such surfaces are dense in the set of Riemann surfaces (\cite{Belyi}).  In contrast to the previous examples , we prove

\begin{Th}\label{thm:1}

Let $S$ be a Belyi surface of genus $g\ge 2$, as $g \to \infty$ 
 the length of the shortest simple closed geodesic on $S$, denoted $\mathrm{syst}(S)$, is bounded by 
\begin{equation}\label{syst}2.809 \le E(\mathrm{syst}(S)) \le 3.085.\end{equation}  
When $E$ is the expected value. In particular, $E(\mathrm{syst}(S))$ is independent of $g$ for surfaces with large genus.

\end{Th}

We also get information about the lengths of some longer geodesics. 
We consider  a set of {\it prime non intersecting} geodesics arranged by length $\gamma_1\le \gamma_2\le \dots \le \gamma_i ,\dots$. We will show that they are all simple and therefore, it is always possible to complete this set and get a pair of pants decomposition \cite{buser}.
We show that 
\begin{Th}\label{thm:2}
For any fixed $i$, $E(l(\gamma_{i}))$ is independent of the genus $g$ as $g \to \infty$.
\end{Th}

This estimate in fact holds for roughly $\frac{g}{\log{g}}$ short geodesics out of the total $3g-3$ which requires to get a pants decomposition of the surface. We provide estimates for the lengths of these geodesics and show that they are all bounded by $ C \log g$. This result is in sharp contrast to Buser's result for general Riemann surfaces where the shortest geodesics are  $\sim \log g$, the same magnitude as  the longest geodesics covered by our estimate.

We analyze the length of the geodesics using the method developed by Brooks and Makover in \cite{BMR}. They use cubic graph to generate the  Belyi surfaces and endow the set  surfaces with probability measure inherit  from random cubic graphs.
 
In section \ref{sec:basconst} we will describe the basic construction of compact surfaces from cubic graph. We review  the connection between the metric structure of the surface and the combinatorics of the graphs.Next, in section \ref{sec:geods} we connect cycles in the graph with geodesics on a surface,
and show that our methods allow us to consider cycles in the graph $\Gamma$ of length $\sim \log{g}$ (see \ref{maxlth}).  We will do so by investigating an interesting explicit example of random matrix multiplication and analyzing the different moments of the distribution of the entries. As it turns out, the entries of the matrix product are a Stern sequence which is an integer sequence of independent interest.

Next, in section \ref{sec:lths} (Theorem \ref{linear}) we use the above results to show that the length of the geodesic associated to such a cycle is growing linearly with the length of the cycle. A detailed analysis of the growth rate of the Stern sequence shows that for a given cycle length on the graph the length of geodesics on the surface is concentrating around the mean length.  We produce the Stern sequence using random products of $2\times 2$ matrices. 

These random products are similar to processes described   by Viswanath (\cite{Visw}) for random Fibonacci sequence, Lima and Rahibe \cite{LR} in the physics of disordered systems, and others. In our example, due to the properties of the Stern sequence, we get more detailed results on the random product then we could find in the literature for similar processes. This example may be of an independent interest in the study of random products of matrices, hence we include a more thorough discussion of our investigations than is necessary to prove our main results.

Finally, in section \ref{sec:short} we calculate the numerical bound to length of the shortest geodesic.

\begin{figure}[htbp] 
   \centering
   \includegraphics[width=5in]{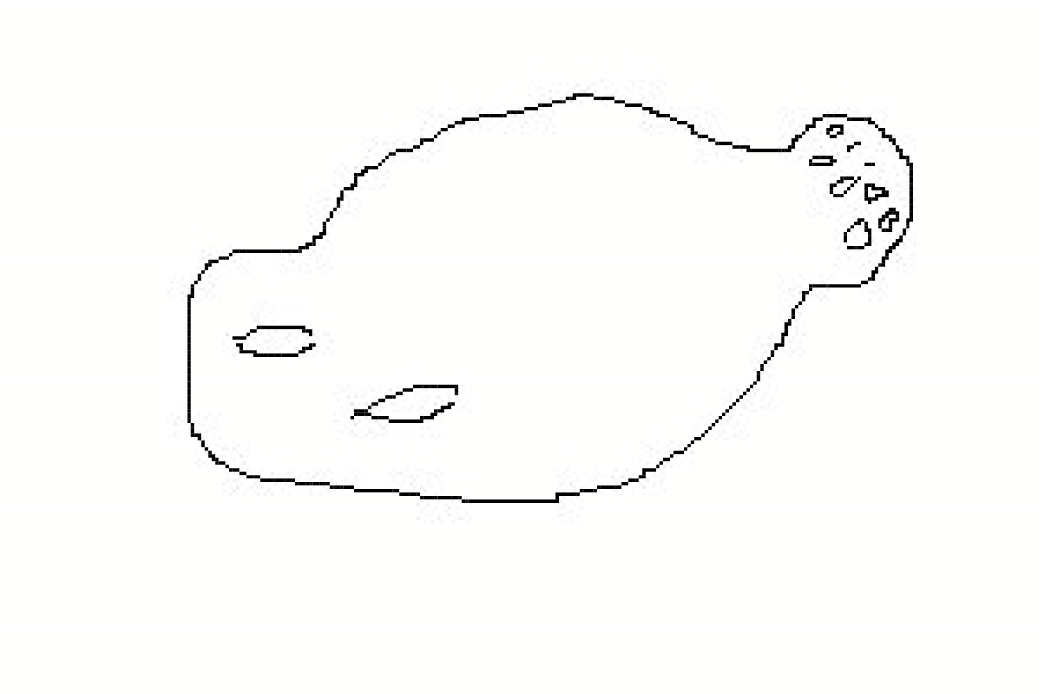} 
   \caption{A ``generic'' Riemann surface.}
   \label{fig:generic}
\end{figure}

This gives a somewhat surprising picture of such random Riemann surfaces.  In (\cite{BMR}) Brooks and Makover showed that the Cheeger constant of  random Riemann surfaces is bounded from below. Hence the geodesics that we find can not disconnect large pieces  of the surface. The picture becomes even more complicated by the result (\cite{BMR,G}) that shows that there is an embedded ball on the surfaces  with area of $ \sim \frac{2}{3}$ of the area of the surface.  A ``generic'' Riemann surface therefore looks something like figure \ref{fig:generic}.

\section{Basic Construction}\label{sec:basconst}
\begin{figure}[htbp] 
   \centering
  \includegraphics[width=2.75in]{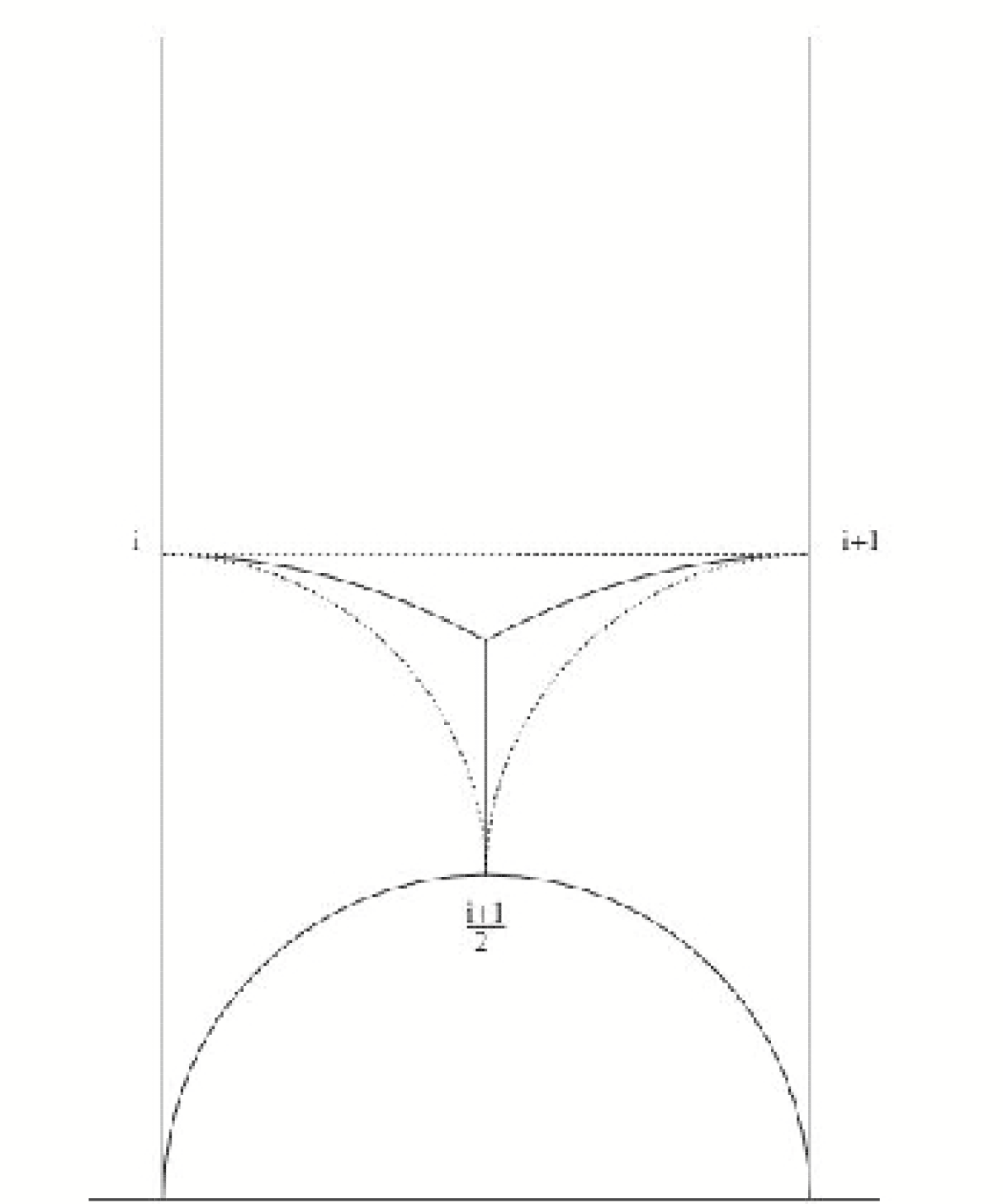}
   \caption{An ideal triangle.}
   \label{fig:idealtri}
\end{figure}

We wish to construct Riemann surfaces from cubic (3 regular) graphs\cite{BMR} .  We will need to add additional combinatorial structure an orientation. We define  orientation $\calO$ on the graph,  as a cyclic permutation of the edges emanating from a each vertex $v$.  
Each vertex in the graph will be identified with an ideal hyperbolic triangle $T$ with vertices at $0$, $1$, and $\infty$ (see Figure \ref{fig:idealtri}).
We may think of the points in $P=\left\{ i,i+1,\frac{i+1}{2}\right\}$ as ``midpoints'' of the corresponding sides of the triangle $T$.  The solid lines in Figure  \ref{fig:idealtri} are geodesics joining the points in $P$ with the point $\frac{1+i\sqrt{3}}{2}$, while the dotted lines are horocycles joining pairs in $S$. Now we can glue neighboring triangles subject to two conditions the first is that the ``midpoints'' of sides are glued together, and the second  is that the orientation $\calO$ at the vertex agree with the orientation of $T$.

%
  Such a gluing is uniquely determined given $(\Gamma,\calO)$.  The resulting surface is an open complete Riemann surface $S^{O}(\Gamma,\calO)$, and we define the compact surface $S^C(\Gamma,\calO)$ as its conformal compactification. 
  While for any given graph of size $n$ there are only $2^{n}$ such surfaces, the following theorem will allow us to model ``generic'' surfaces using surfaces generated from graphs.
\begin{Th}[\cite{Belyi}] The set of surfaces constructed from cubic graphs is dense in the set of Riemann surfaces.\end{Th}  

To work with the set of surfaces generated from random cubic graphs with $n$ vertices, we will need a probability measure for the set of such graphs.  While the fact that for each $n$ we have a finite set of such graphs,  it is hard to work with \cite{Wormald} \cite{bollobas} therefore we will use Bollobas' configuration model.

Put $6n$ balls in a hat; label the balls using the numbers $1,2, \dots ,(2n)$, with three copies of each number. Then pick at random pairs of balls.  We can define a graph by taking a set of $1,2, \dots ,(2n)$ vertices, and connecting vertices  $v_i$ to $v_j$ by an edge if a pair of two balls marked with $i$ and $j$ have been picked together.  We will endow the set of oriented cubic graphs with a probability measure by picking a graph using the configuration model and then flipping an unbiased coin at each vertex to pick an orientation. Note, we will allow loops and double edges since they not do not interfere with the construction of the surface $S^O(\G, \calO)$. We use the notation $\calFn$ for the set of cubic graphs on $2n$ vertices with the above probability measure, and $\calFns$ for the set of
 {\emph oriented} cubic graphs with the same probability measure.  

For the unoriented graphs Bollobas proved (\cite{bollobas} \cite{MWW})
\begin{Th}\label{poisson} Let $X_i$ denote the number of
closed paths in $\Gamma$ of length $i$. Then the random variables $X_i$ on
$\calF_n$ are asymptotically independent Poisson distributions with means
$$\lambda_i = {{2^i}\over{2i}}.$$

\end{Th}
%

In order to study the lengths of simple closed geodesics on the compact surfaces $S^{C}(\Gamma,\calO)$, we must understand the relation of the metric structure of $S^{C}(\Gamma,\calO)$ and $S^{O}(\Gamma,\calO)$.  The following two theorems show that as $n \to \infty$, almost all surfaces $S^{C}(\Gamma,\calO)$ have a global  metric structure arbitrarily close to $S^{O}(\Gamma,\calO)$.
\begin{Th}[{\cite{platonic}}]\label{PS}
 We say that $S^O(\G,\calO)$ has cusps $\ge L$ if    there is a  set of disjoint  horocycles, one around each cusp and each has length $\ge L$. Then for every $\eps$, there exists numbers $L, r,$ 
and $y$ such
that, if the cusps of $S^O$ have length $\ge L$, outside the
union of cusp neighborhoods $\U= \cup_{i=1}^k f_i^{-1}(C^{y})\subset 
S^{O}$ of 
the cusps $C^{i}$, and $\V= \cup_{i=1}^k B_{r}(p_i)\subset S^{C}$, the metrics
$ds_C^2$ and $ds_O^2$ satisfy
$${{1}\over{(1+ \eps)}} ds_O^2 \le ds_C^2 \le (1 + \eps)ds^2_O.$$
\end{Th}

The large cusp condition is necessary and enables us to compare the global metric of open and compact surfaces.
Let $Q$ be a property of $3$-regular graphs with orientation, and denote by
$\Prob_n[Q]$ the probability that a pair $(\Gamma, \calO)$ picked from
$\calFns$ has property $Q$.

\begin{Th} [\cite{BMR}]\label{cusps}

For every $L>0$, as $n \to \infty$, 
$$\Prob_n [ S^O(\Gamma, \calO)\ {\hbox{has cusps of length}} >L\ ] \to 1.$$

\end{Th}

\section{Geodesics on $S^{C}(\Gamma,\calO)$.}\label{sec:geods}

The discussion in the previous section shows that we can use the graph $\Gamma$ to get information about the surface $S^{C}(\Gamma,\calO)$.  To do this, we must begin with a cycle in $\Gamma$ and its' associated geodesic on $S^O(\G, \calO)$.  Next, we track the geodesics as the cusps of $S^{O}(\Gamma,\calO)$ are closed to give $S^{C}(\Gamma,\calO)$. 

 The  geodesics of  $S^O(\G, \calO)$ are described in the oriented graph $(\G, \calO)$ as follows; let $\calL$ and $\calR$ denote the
matrices 
$$\calL = \btbt 1&1\\0&1 \etbt \quad \calR = \btbt 1 & 0 \\ 1 & 1
\etbt.$$
A closed path $\calP$ of length $k$ on the graph may be described by
starting at the midpoint of an edge, and then giving a sequence $(w_1,
\ldots, w_k),$ where each $w_i$ is either $l$ or $r$, signifying a left
or right turn at the upcoming vertex. We then consider the matrix
\begin{equation}\label{MP}M_{\calP}(k) = W_1 \ldots W_k,\end{equation}
where $W_j = \calL$ if $w_j = l$ and $W_j = \calR$ if $w_j=r$. 

The closed path $\calP$ on $\Gamma$ is then homotopic to a closed
geodesic $\gamma(\calP)$ on $S^O(\Gamma, \calO)$ whose length
$\l(\gamma(\calP))$ is given by 
\begin{equation}\label{length}2\cosh({{\l(\gamma(\calP))}\over{2}})= \tr(M_{\calP}).\end{equation}

We have to check what happen to the geodesics on $S^O(\G,\calO)$  as we close the cusps and get $S^O(\G,\calO)$.
 Given a simple geodesic $\gamma$ on $S^O(\G,\calO)$ let $\breve{\gamma}$ be the image of $\gamma$ in $S^C(\G,\calO)$. 
  
 First $\breve{\gamma}$ might be homotopicaly trivial. 
This can happen in two ways.  One is that the orientation on the path is uniform (all "Left" or all "Right") and in this case the path is circling a cups and the length formula (\ref{length})  gives $0$ since there is no geodesic in its homotopy class in $S^O(\G,\calO)$. The second possibility is that  $\breve{\gamma}$ is nontrivial on $S^O(\G,\calO)$ but bounds a disk in $S^C(\G, \calO)$, in this case the corresponding cycle on the graph disconnects the graph.

Second, the image of two non equivalent geodesics on $S^C(\G,\calO)$ might become homotopicaly equivalent on  $S^C(\G,\calO)$. The treatment of this case is similar to the second case above since,  in this case the images of the two geodesics bound a cylinder in $S^C(\G, \calO)$ and therefore the cycles on the graph disconnect the graph.

The probability that one or two cycles will disconnect the graph tends to $0$  by the following
\begin{Lm}\label{vanish}
Let $\calP$ be a cycle or  union of two cycles  on $\G_n$ with $\length(\calP) < C\log n$ then

$$\Prob_n [ \calP\ {\hbox{disconnect $\G_n$}} \ ] \to 0.$$
\end{Lm}
The proof of this lemma is a straight forward application of the following results

\begin{Th}[\cite{BOL}]
$$\Prob_n [h(\Gamma) >\frac{2}{11}] \to 1\ {\hbox{as}}\ n \to \infty.$$
\end{Th}
Therefore if $\calP$ disconnects $\G$ it induces a subgraph $H \subset \G$ with $v(H) \le Const*v(\calP)$
but it is known that
\begin{Th}[\cite{Wormald,  bollobas}]
Let $H$ be a graph such that $v(H)<e(H)$.  Then the expected number of copies of $H$  in a random cubic graph with $n$ vertices is $O(n^{-1})$
\end{Th}

Given  a geodesic $\gamma$ on $S^O(\G, \calO)$. Let $\breve{\gamma}$ be the image of $\gamma$ in $S^C(\G, \calO)$ then if  $\breve{\gamma}$ is not homotopicly trivial after closing all cusps, then there is a unique geodesic $\hat{\gamma}$ homotypical to $\breve{\gamma}$ in $S^C(\G,\calO)$. Under  the large cusps condition we can compare the length of $\gamma$ and $\hat{\gamma}$
\begin{Th}[\cite{platonic}]

$$\frac{\l(\gamma)}{1+\epsilon} \le \l(\hat{\gamma}) \le \l(\gamma)$$
\end{Th} 

\begin{figure}[htbp] 
   \centering
   \includegraphics[width=3.5in]{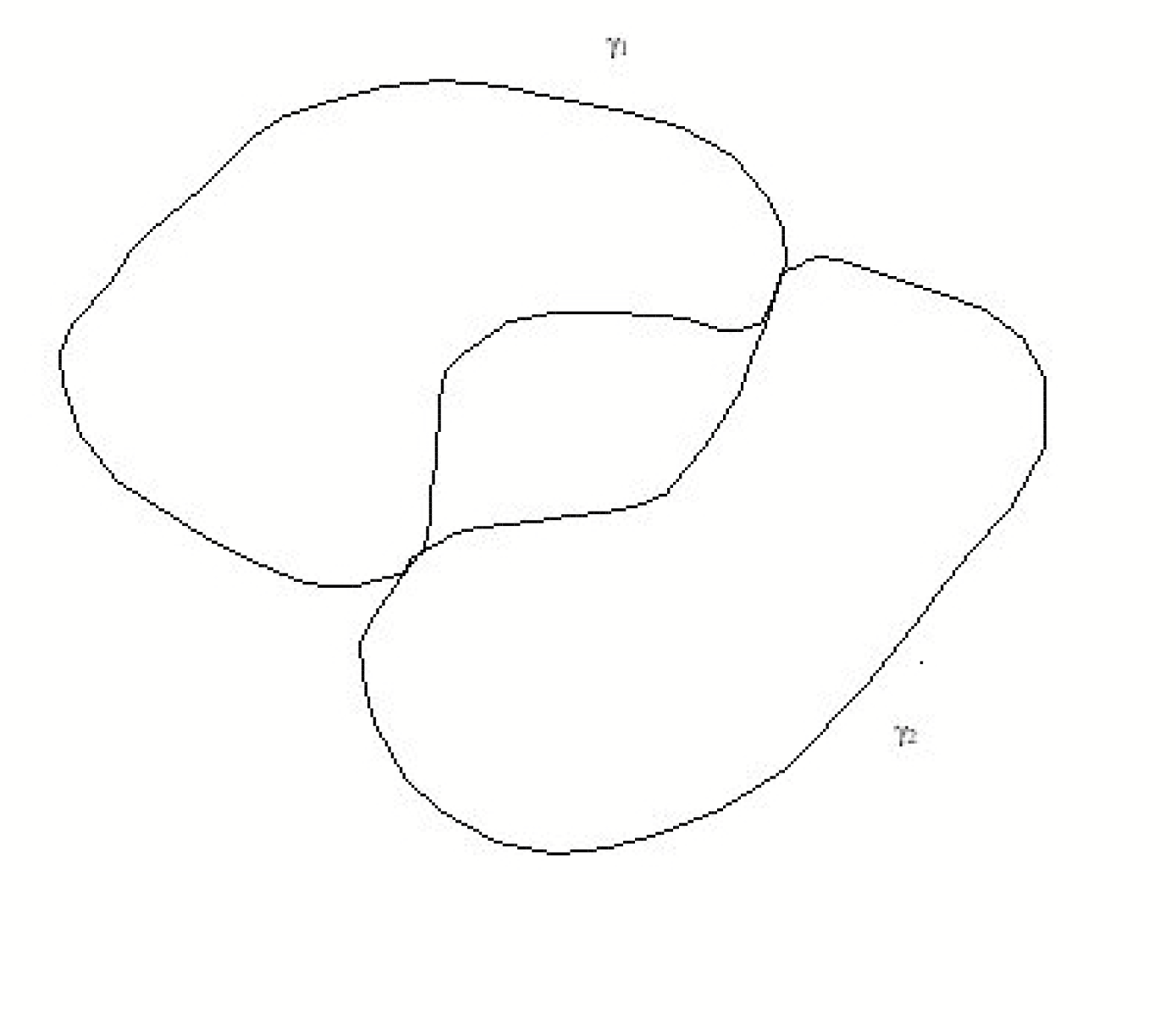} 
   \caption{We cannot allow two intersections.}
   \label{fig:intgeos}
\end{figure}

Before computing the expected length of {\bf the} geodesic corresponding to a cycle of length $l$ on the graph we must deal with the question of how far the Bollobas estimate of Poisson distributions holds. It is clear that for very long cycles the distribution is different \cite{Garmo}.  The distribution of Hamiltonian cycles is known \cite{Garmo} . At the other extreme, for very short cycles the Poisson distribution is a very good estimate. Recently McKay, Wormald, Wysocka  examine this problem \cite{MWW} and gave a very detailed estimate for a more general problem. We will modify part of their results to  suit our needs. In particular, their result covers the distribution of non-intersecting cycles. We will allow a very limited intersection, namely we allow for two cycles to touch once. In this case we can easily ``pull" the corresponding path apart on the surface and thus the resulting geodesic on $S^C(\Gamma, \calO)$ will be non-itersecting. Note that allowing more then one component in the intersection of two cycles will result in counting some geodesics more then once, as seen in \ref{fig:intgeos} .


We use the following to determine how long cycles may be so that the probability that the intersection of any two cycles will have two or more components goes to 0 as $n \to \infty$.
\begin{Th}\label{thm:cycint}
Let $K$ be a cubic graph, and $\mathcal{C_{r}}$ be the collection of all cycles of length $r$ in $K$, with $3 \le r \le \alpha$.  Let $\mathcal{C} \subset \mathcal{C_{r}} \times \mathcal{C_{s}}$, $s \le r$, be all pairs of cycles $(C_{1},C_{2})$ with $C_{1}\bigcap C_{2}\neq \emptyset$ and $C_{1}\neq C_{2}$.  In addition, suppose that the number of components $p$ in $C_{1}\cap C_{2}$ is at least 2.  Given a cubic graph $\Gamma$ with $n$ vertices, \begin{equation}\label{intprob}\sum_{(C_{1},C_{2})\in \mathcal{C}}\mathbf{P}(C_{1}\bigcup C_{2} \subseteq \Gamma) \le O(1)\sum_{j\ge 1, p\ge2}\frac{(2\alpha^{3})^{p-1}}{(p-1)!^{2}}n^{r+s-p-j}\left(\frac{2}{n}\right)^{r+s-j}\end{equation} where
 $j$ is the number of edges in $C_{1}\cap C_{2}$.  \begin{proof} Follows directly from formula 2.7 in (\cite{MWW}).\end{proof}\end{Th}
Summing for all pairs $r,s \le \alpha$, and assuming that $\alpha$ is chosen so that $\alpha^{3}=o(n)$,  we find that the probability that any cycle of length less than $\alpha$ has intersection with more than one component with a different cycle also of length less than $\alpha$ is \begin{equation}\label{alphadef}O\left(\frac{2^{2\alpha-1}}{n^{2}}\right)\end{equation}  If we choose $\alpha$ so that $$2^{2\alpha -1} = o(n^{2})$$ then this probability will go to zero as $n \to \infty$.  Thus, \begin{eqnarray}2^{2\alpha -1} &=& o(n^{2}) \implies \\ 2^{2\alpha -1} &=&n^{2-\epsilon}\\2\alpha - 1 &=&(2-\epsilon)\log_{2}{n}\\\alpha &=& (1-\epsilon')\log_{2}{n}\label{maxlth}\end{eqnarray}

\section{Lengths of Geodesics}\label{sec:lths}

We consider cycles of length $N$ in the graph, and the geodesics on the surface associated with each cycle.  The length of the geodesic associated to any given cycle can be computed using (\ref{length}).  Clearly these lengths will depend upon the orientation of the vertices, and for cycles of length $N$ there are $2^{N}$ possible orientations.  We need to understand the distribution of the lengths of the associated geodesics.

We begin by considering the expected value and the distribution for the traces of the matrices $M_\mathcal{P}$ from (\ref{MP}), where the paths $\mathcal{P}$ are some random set of paths of length $N$.  We will do this by considering the individual matrix entries.
Let $M_{\mathcal{P}}(i) = W_{1}W_{2}...W_{i}$ with $i \le N$, and $M_{\mathcal{P}}(0) = Id$, the $2 \times 2$ identity matrix.  If $$M_{\mathcal{P}}(i) = \begin{pmatrix}
    a  &   b \\
  c    &  d
\end{pmatrix},$$ then $M_{\mathcal{P}}(i+1)$ is either $$\begin{pmatrix}
    a+b  &   b \\
  c+d    &  c
\end{pmatrix}$$ or $$\begin{pmatrix}
    a  &  a+ b \\
  c    &  c+d
\end{pmatrix}.$$  Consider the top row of $M_{\mathcal{P}}(i)$.  At step $i=0$ it is $(1\,0)$, and thus 1 and 0 are the only values possible for either entry.  At step $i=1$, the top row is either $(1\,0)$ or $(1\,1)$, and as $i$ continues to increase, it is clear that the possible values are values from the previous step, or sums of neighboring values from the previous step.  Thus, we can build a table of possible values for matrix elements by starting with 1 and 0, putting the sum 1 in between, then putting the sums 2 and 1 in the new spaces, \ldots 
\begin{equation}\label{sternseq}
1\cdots 2 \cdots 1 \cdots 1 \cdots 0
\end{equation}

After $i$ steps we have a sequence with $2^{i}+1$ entries.  This sequence is one example of a {\it Stern sequence} (\cite{Stern,GG1,GG2,GG3}).  Stern sequences have many nice properties, and the moments at any step $i$ can be computed directly.  For example, each term, except for the initial terms $a$ and $b$, shows up as part of two new terms in the next sequence, so if $S(i)$ is the sum of terms at step $i$, then $S(i+1) =3 S(i) - (a+b)$.  
Stern sequences, in the guise of the {\it Stern-Brocot tree} also show up in the work of Viswanath (\cite{Visw}) on random Fibonacci sequences, where he considered random products of the matrices $$A = \btbt 0&1\\1&1 \etbt \quad B = \btbt 0 & 1 \\ 1 & -1
\etbt.$$  

We need only consider values of diagonal elements of $M_{\mathcal{P}}(i)$.  At any given step, one diagonal element will change and the other will not, so it is clear that we need to understand both the distribution of the diagonal elements and the dependence between them to determine the distribution of values of the trace.  The rows of the matrix associated with  a path $\mathcal{P}$ will be  generated by Stern sequences starting with $(1\,0)$ and $(0\,1)$, which we denote $S_{u}$ and $S_{l}$.  .  These are simply reflections of each other, and we omit the leading (or trailing) zero from the sequence.  Counting from $i=0$, we see that the sum of elements is $$S(i)=\sum_{j=1}^{i-1}3^{i}+2,$$ while the number of elements $N(i)=2^{i}$, hence the expected value at step $i$ is \begin{equation}
\label{sternmean}
E_{i}=\left(\frac{\sum_{j=1}^{i-1}3^{i}+2}{2^{i}}\right)=\frac{3^{i}+1}{2^{i+1}}
\end{equation}

We can compute the variance $\sigma^{2}_{i}$ by determining the sum of the squares of the elements in the sequence, $S_{2}(i)$, as follows (\cite{email}).  Consider three neighboring terms in the sequence $$x_{j}+x_{j+1}\qquad x_{j+1}\qquad x_{j+1}+x_{j+2}.$$  Setting $$A(i) = \sum_{j=1}^{2^{i}} x_{j}^2$$ and $$
B(i) = \sum_{j=1}^{2^{i}} x_{j}x_{j+1}$$ one gets the recurrence relations \begin{eqnarray}
A(i+1) &=& 3A(i) + 2B(i)-2\\ B(i+1)&= &2A(i) + 2B(i)-2\\\implies A(i) &=& 5A(i-1) - 2A(i-2)-1
\end{eqnarray}

A messy but straightforward computation gives \begin{equation}
\label{ }
A(i)=
{\textstyle \frac{2^{-x-1} \left(17 2^x
   \left(-5+\sqrt{17}\right)+\left(5+\sqrt{17}\right)^x
   \left(-34+6 \sqrt{17}\right)+\left(5-\sqrt{17}\right)^x
   \left(-51+11 \sqrt{17}\right)\right)}{17
   \left(-5+\sqrt{17}\right)}}\nonumber
\end{equation} using $A[1]=2$, $A[2]=7$.
Similar calculations give the sums of $n$th powers of terms of the sequence.  Thus we can compute any central moments of our Stern sequence, and the distribution of diagonal matrix entries.

We must convert our knowledge of the moments of the diagonal entries into information about their sum.  Since means are additive, we can compute the mean value of the trace using (\ref{sternmean}),\begin{equation}
\label{tracemean}
E_{Tr,i}=2\left(\frac{\sum_{j=1}^{i-1}3^{i}+2}{2^{i+1}}\right)=\frac{3^{i}+1}{2^{i}}
\end{equation}

 Clearly a dependence exists between the value of the upper diagonal and lower diagonal elements $x_{u}$ and $x_{l}$.  Hence, determining the variance of the traces will require a computation of the covariance of the Stern Sequences for the upper and lower diagonal elements.  The covariance can be computed using \begin{equation}
\label{covdefn}
Cov(S_{l},S_{u}) = <s_{l}s_{u}> - <S_{l}><S_{u}>
\end{equation} where $<>$ indicates the expected value.  Equation (\ref{sternmean}) gives us the second term, so we need to compute the mean of products of associated elements in the two sequences, in other words the product of the diagonal entries for each matrix $M_{\mathcal{P}}$ corresponding to paths of length $i$.

\begin{figure}[htb]
\includegraphics[width=5in]{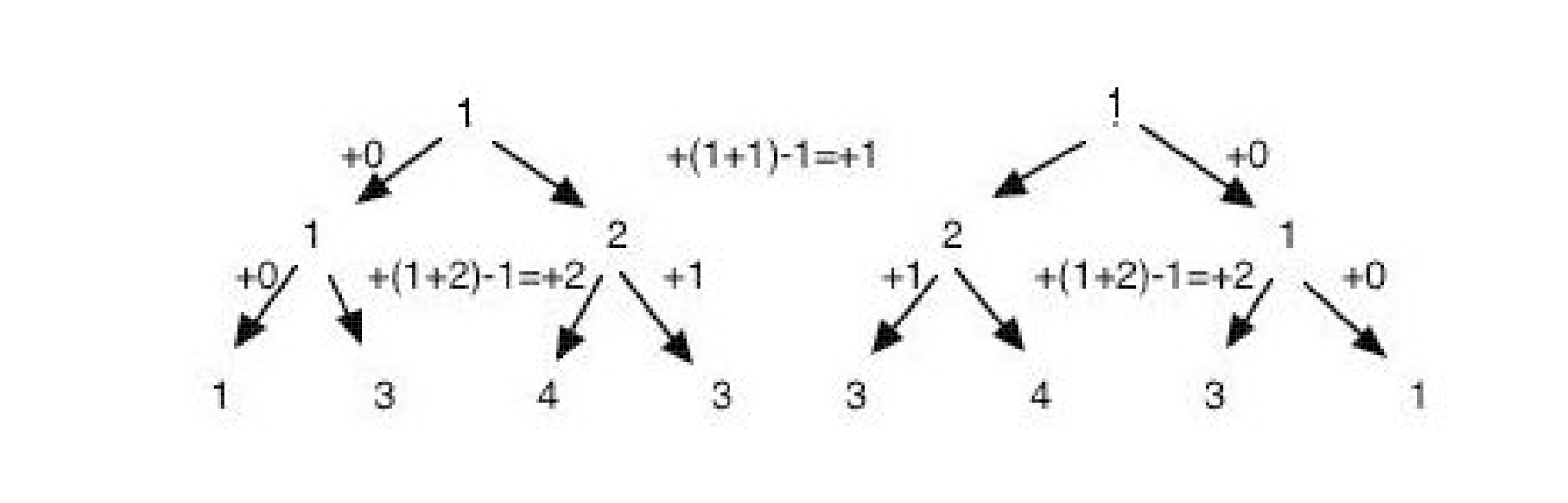}
\caption{Computing the products of corresponding elements of the Stern sequence.}
\label{mult}
\end{figure}

We compute the product of diagonal matrix elements using the diagram in Figure \ref{mult}, which follows directly from an investigation of the matrices.  After an initial choice of either $\calL$ or $\calR$, the diagonal product is 1.  Each number in a given row will generate two children, corresponding to a choice of $\calL$ or $\calR$. We determine the value of the children using two simple rules.  If the path to the child is a continuation of a previous path from higher in the diagram, we add the same value that was previously added, otherwise we add the parents value to its' siblings, subtract one, and that becomes the addend for the new direction.  The outer edges of the tree have addends of zero, which correspond to paths containing only $\calL$ and only $\calR$.

\begin{figure}[htb]
\begin{center}
\includegraphics[width=4in]{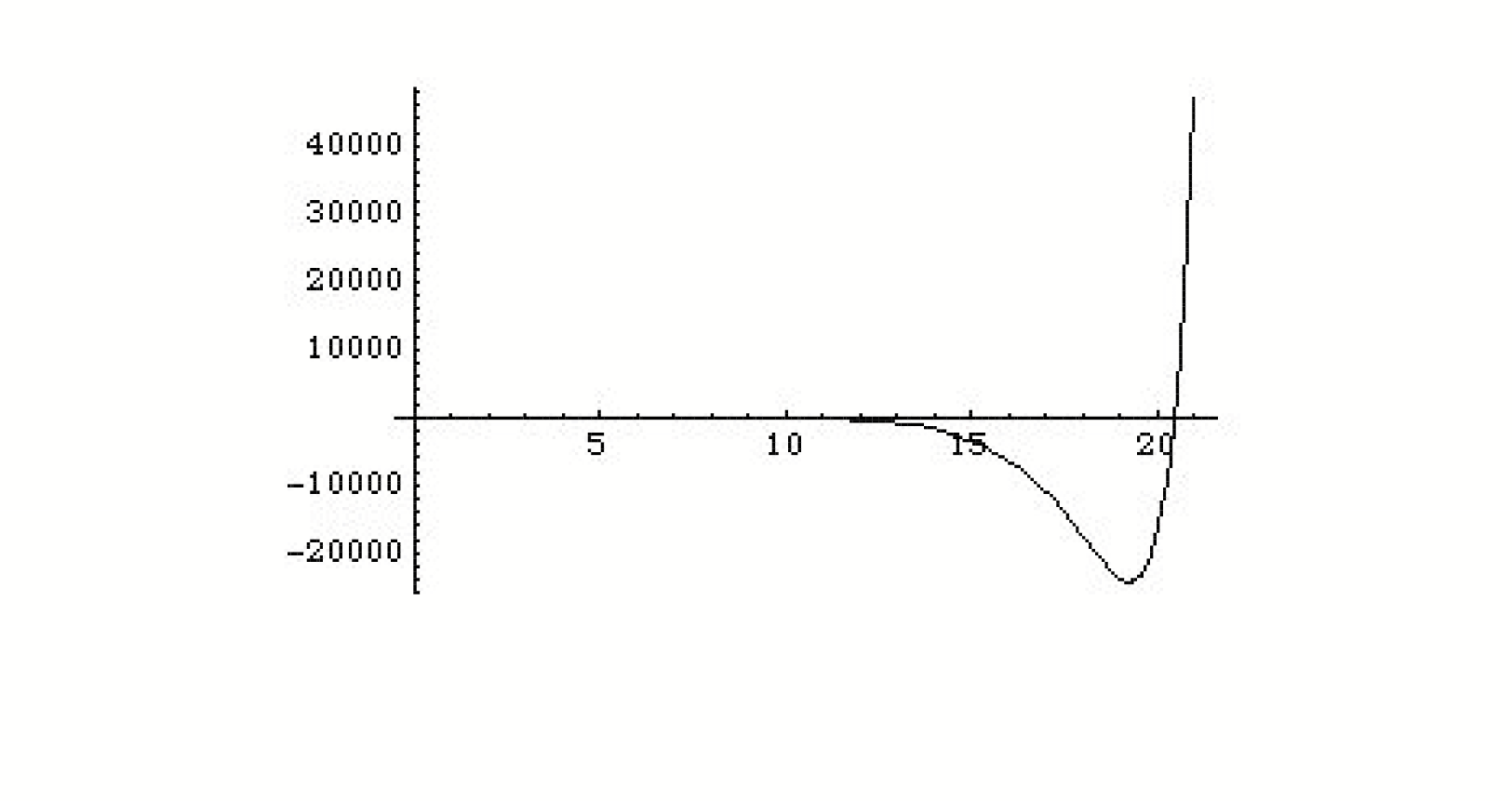}
\caption{ The covariance for the two diagonal elements of a matrix $M_{\mathcal{P}}(n)$.}
\label{covplot}
\end{center}
\end{figure}

As in the computation for the sums of squares, we get a recurrence relation  \begin{equation}
C(i+1) = 5C(i) + 2C(i)-2^{i-1}
\end{equation}
where $$C(i) = \sum_{i=1}^{2^{i}} x_{u,i}x_{l,i}$$ and $C(1)=2$, $C(2)=6$.  Thus, \begin{equation}
\label{prodformula}
C(x)={\scriptstyle \frac{1}{17} 2^{-x-2} \left(17 2^{2
   x+1}-\left(-17+\sqrt{17}\right)
   \left(5+\sqrt{17}\right)^x+\left(5-\sqrt{17}\right)^x
   \left(17+\sqrt{17}\right)\right)}
\end{equation}  Substituting equations (\ref{sternmean}) and (\ref{prodformula}) into (\ref{covdefn}) gives \begin{equation}
\label{covariance}
Covariance={\textstyle \frac{17 \left(-1+2^{2 x+1}-2
   (3^x)-9^x\right)- \left(-17+\sqrt{17}\right)
   \left(5+\sqrt{17}\right)^x+\left(5-\sqrt{17}\right)^x
   \left(17+\sqrt{17}\right)}{17 (4^{x+1})}}
\end{equation}Looking at a plot of (\ref{covariance}) in Figure \ref{covplot}, we see that while for short cycles there is a negative correlation between the diagonal elements, this quickly changes over to a positive one.  In fact, we can compute the correlation between the two diagonal matrix entries, and see that $$\lim_{n\to \infty}Cor = \frac{51-11 \sqrt{17}}{34-6 \sqrt{17}} \simeq .61$$

Knowledge of the covariance for the upper and lower diagonal elements of the matrices allows us to compute the variance for the traces, which gives us a basic understanding of the distribution of traces for cycles of length $N$.  
\begin{equation}
\label{trvar }
Variance = {\scriptstyle 4^{-x} \left(-1+2^x-2
   3^x+4^x-9^x+\left(5-\sqrt{17}\right)^x+\left(5+\sqrt{17
   }\right)^x\right)}
\end{equation}
Unfortunately, the lengths of geodesics on the surface corresponding to the given cycles is determined using (\ref{length}), \begin{equation}
\label{length2}
Length(\gamma)= 2\cosh^{{-1}}(\frac{Tr}{2})=2\ln{\left(\frac{Tr+\sqrt{Tr^{2}-1}}{2}\right)}.
\end{equation}
When one has two independent random variables related by some function $Y=f(X)$, the standard technique used to deduce information about the distribution of the $Y$'s given information about the distribution of the $X$'s is to use a Taylor approximation, and evaluate at the expected value of $X$, $E(X)$.\begin{eqnarray}
Y & \approx &f(c)+f'(c)(X-c)+\frac{f''(c)(X-c)^{2}}{2}\qquad {\text Let }c=E(X)\implies\\
Y& = & f(E(X)) + f'(c)(X-E(X)) +  \frac{f''(E(X))(X-E(X))^{2}}{2}.\label{expY}
\end{eqnarray} One can now approximate $E(Y)$ using (\ref{expY}).  The second term disappears, and since in the case we are considering $f(X) \approx \ln(x)$, the third term gives roughly $$-\frac{Var(X)}{E(X)^{2}}$$ which approaches $-\infty$ exponentially.  Thus, we get no lower bound even for the mean of the lengths as $n$ increases.

The following Lemma, which follows directly from the arithmetic geometric mean inequality, shows that the means of the lengths cannot be growing any faster than linearly in $n$.   
\begin{Lm}\label{upper}
For a positive  random variable $X$  $$E(\log(X)) \le \log(EX)$$
\end{Lm} 

We need to give a lower bound for the growth of the lengths.  The following proposition shows that the lower bound is also linear.
\begin{Th}\label{lowerbound}
The proportion of the standard Stern sequence which is growing exponentially in $N$ as $N \to \infty$ approaches 1.
\end{Th}
\begin{proof}
First, recall that any pair of neighbors in the $N$th Stern sequence represents a row in $\Mp(N)$.  Next, note that in $\Mp(1)$ one row must be $(1,1)$, without loss of generality we will assume that this is the bottom row.  Thus, at the expense of one step, we can consider the growth of the standard Stern sequence beginning with $1 \qquad 1$, which contains $2^{k}$ pairs after $k$ steps.

We may now consider a Bernoulli process defined as follows: for any pair in the $k$th Stern sequence, consider the four pairs which are its' children in the $k+2$nd Stern sequence.  These come from the four choices of turns possible, $\calR \calR$, $\calR \calL$, $\calL \calR$, and $\calL \calL$.  Since we are considering the trace of $\Mp$, we need only consider the right element in each pair.  If our initial pair is $a \qquad b$, we must consider two possibilities, $a>b$ or $a<b$.  The possibility that $a=b$ will only occur in the last pair, which has vanishing probability. 

In Figure (\ref{CSstring}), we see that the four possible values for the diagonal element are $2a+b$, $a+b$, $a+2b$, and $b$.  If $a < b$, then one of the four is double the original diagonal element $b$, while if $a > b$, three of the four are double the original.  As $N \to \infty$, the probability that $a<b$ and $a>b$ approach 50\%.  Thus, defining success for our Bernoulli process as a doubling, the probability of success is $1/2 - \epsilon(N)$, with $\lim_{N\to \infty}{\epsilon}=0$.
\begin{figure}
\begin{center}
\includegraphics[width=4in]{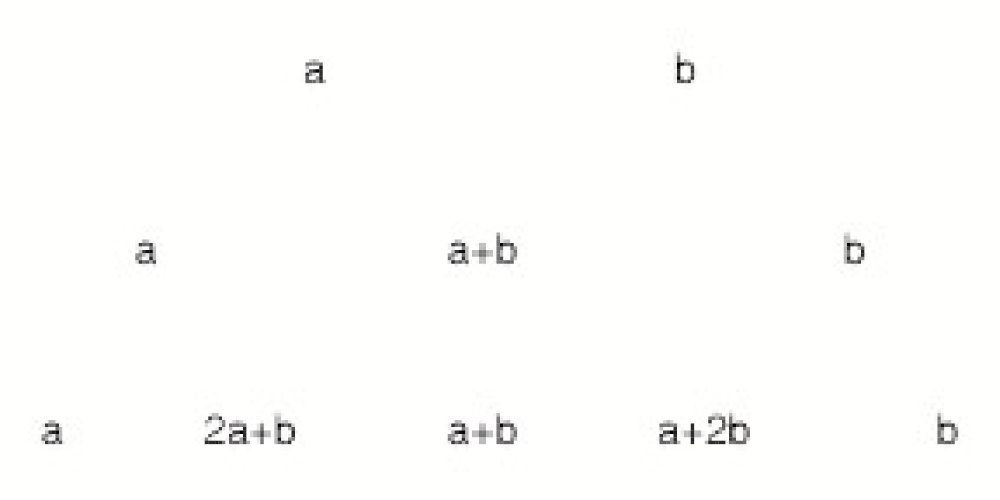}
\caption{Three out of four entries at least double every two steps}
\label{CSstring}
\end{center}
\end{figure}

Considering paths of length $N=2k+1$, the Bernoulli process will have $k$ steps.  As $N \to \infty$, the  mean will approach $N/2$ and the probability that there will be at least $k$ successes will approach one.  Therefore the probability that any entry in the $n$th Stern sequence is at least $2^{k}$ approaches 1 as $N \to \infty$.\end{proof}

Keeping better track of the minimal growth using a table like Figure \ref{CSstring} allows us to give a numerical estimate for the lower bound on the growth rate.  For example, if $b<a$, then of four choices for the second term in a pair of Stern sequence entries, two have at least tripled ($2a+b$ and $b+2a$).  One then has a multinomial process, and if the number of steps corresponding to a single trial is allowed to increase, the lower bound also increases.

Allowing five steps, one calculates that a lower bound for the growth factor of almost every element is given by $$2^{37/320}3^{1/16}5^{3/80}7^{1/40}11^{1/80}13^{1/160} \approx 1.35502.$$  Such calculations quickly get too long to do by hand, and the return on investment becomes minimal since the mean is growing like $1.5^{N}$.  A simple Mathematica routine can push the calculations significantly farther, however; for example if the basic step size is 15, then the growth factor is approximately $1.43925$.

Combining this with (\ref{length2}) and Lemma \ref{upper}, we get \begin{Th}\label{linear}As $N \to \infty$, the expected value of the length of the simple closed geodesics associated to cycles of length $N$, $E_{N}(l)$, is bounded by $$(\log{1.43925})N \le E_{N}(l) \le (\log{1.5})N.$$\end{Th}

 \section{The Length of the Shortest Geodesic}\label{sec:short}

 In this section we will use  the results from previous section to estimate
the  length of the shortest closed geodesic on $S^C(\G, \calO)$.  The length
of the  shortest  closed geodesic is an important geometric invariant since
it is twice  the  injectivity radius. Let $\syst(S^C)(\G, \calO)$ be length
of the shortest closed geodesic on $S^C(\G, \calO)$, we will show that :

 \begin{Th}
$$2.809\le  \syst(S^C)(\G, \calO) \le 3.085$$ \end{Th} 

It is interesting to compare this result in \cite{BMR} with the the
behavior  of the  largest embedded ball on$(S^C)(\G, \calO)$  which gives a
linear growth to the largest embedded ball.

We will start with the  upper bound.
As we have seen on a random graph the distribution of short cycles is
independent on the size of the graph and it is Poisson distribution with
mean $ \lambda=\frac{2^k}{2k}$ where $k$ is the length of the cycle.
Therefore the probability of not having a $k$ cycle is $e^{-\frac{2^k}{2k}}$
and the probability of having at least one $k$-cycle on the graph is
$(1-e^{-\frac{2^k}{2k}})$. 

As we transition from the graph $\G$ to the compact
surface $S^C(\G, \calO)$, by Lemma \ref{vanish} there are only 2  orientations
that will make the curve that corresponds to a cycle null homotopic.  Hence
the probability that there is a  $k$-cycle  that gives rise to non
trivial geodesic on $S^C(\G,\calO)$ is
$p(k)=\frac{2^{k - 2} - 1}{2^{k - 2}}(1-e^{-\frac{2^k}{2k}})$ , and the
probability that there is no $k$-cycle  that  produces a geodesic on $S^C(\G,
\calO)$ will be $(1-p(k))$. The expected value for the trace  of a $k$ cycle
is $tr(k)\frac{3^k+1}{2^k}$ and the expected value for the  length of the
geodesic  is bounded from above by  by  $2\cosh^{-1}(\frac{tr(k)}{2})$.
$$E( \syst(S^C)(\G, \calO)) \le \sum_{k=3}^\infty
(p(k)\prod_{j=2}^{k-1}(1-p(j)))2\cosh^{-1}(\frac{3^k+1}{2^{k+1}})$$

$$=\sum_{k=2}^\infty  \bigg (\frac{2^{k - 2} - 1}{2^{k -
2}}(1-e^{-\frac{2^k}{2k}})\prod_{j=2}^{k-1}\Big(1-\frac{2^{j - 2} - 1}{2^{j
- 2}}(1-e^{-\frac{2^j}{2j}}) \Big)\bigg)2\cosh^{-1}(\frac{3^k+1}{2^{k+1}})$$

It is easy to see that this series converges rapidly  and we can get a
numerical estimate that the value  is $\sim 3.085$

To get a lower bound we can replace $\log E (tr(M))$ by
$E(2\cosh^{-1}(\frac{tr(M)}{2}))$ and use the rapid decay for the probability
that a graph has large girth.   We calculate
$$\sum_{k=2}^20 \bigg (\frac{2^{k - 2} - 1}{2^{k -
2}}(1-e^{-\frac{2^k}{2k}})\prod_{j=2}^{k-1}\Big(1-\frac{2^{j - 2} - 1}{2^{j
- 2}}(1-e^{-\frac{2^j}{2j}}) \Big)\bigg)E(2\cosh^{-1}(\frac{M}{2}))$$
and get the estimate of $2.809$.

\bibliographystyle{amsplain}
	\bibliography{geo}

\end{document}